%% file: cone_prog_pert.tex
\documentclass[11pt]{article}
\usepackage{amsmath}
\usepackage{hyperref}
\usepackage{fullpage}
\usepackage[dvipsnames]{xcolor}
\usepackage[capitalize, nameinlink]{cleveref}
\usepackage{autonum}

\newcommand*{\tran}{^T}
\newcommand{\RR}{{\mbox{\bf R}}}
\DeclareMathOperator{\minimize}{minimize}

\providecommand{\Id}{I}

\newcommand{\sign}{\ensuremath{\operatorname{sign}}}

\input{defs}

\title{Differentiating Through a Cone Program}
\author{Akshay Agrawal \and Shane Barratt \and 
Stephen Boyd \and Enzo Busseti
\and Walaa M.\ Moursi\thanks{Authors listed in
 alphabetical order.}}

\begin{document}
\maketitle

\begin{abstract}
We consider the problem of efficiently computing the derivative of the
solution map of a convex cone program, when it exists.
We do this by implicitly differentiating the residual map for its
homogeneous self-dual embedding, and solving the linear
systems of equations required using an iterative method.
This allows us to efficiently compute the derivative operator, and its
adjoint, evaluated at a vector.
These correspond to computing an approximate new solution, given
a perturbation to the cone program coefficients (\ie, perturbation
analysis), and to computing the gradient of a function of the solution
with respect to the coefficients.
Our method scales to large problems, with numbers of coefficients
in the millions.
We present an open-source Python implementation
of our method that solves a cone program and returns the derivative
and its adjoint as abstract linear maps; our implementation
can be easily integrated into software systems for automatic differentiation.
\end{abstract}

\section{Introduction}
\label{sec:introduction}

A \emph{cone program} is an optimization problem in which the objective is to
minimize a linear function over the intersection of a subspace and a convex cone.
Cone programs include linear programs, second-order cone programs, and
semidefinite programs.
Indeed, every convex optimization problem can be expressed as a cone
program \cite{nemirovski2007advances}.

Cone programs can be efficiently solved by a number of
methods, including interior-point methods \cite{nesterov1994interior}
and operator splitting methods such as the alternating directions method
of multipliers \cite{ocpb:16}.
Cone programs have found applications in many
areas, including control \cite{boyd1994linear}, machine learning
\cite{hastie2009stats, boyd2011distributed}, finance
\cite{markowitz1952portfolio, boyd2017multi}, supply chain management
\cite{bertsimas2004robust,ben2005retailer}, and energy management
\cite{moehle2018dynamicnot}, to name just a few.
Cone programs are widely used in conjunction with
several popular software systems for convex optimization,
which reformulate a convex optimization problem expressed in a
domain specific language into an equivalent cone program. These include YALMIP
\cite{lofberg2004yalmip}, CVX \cite{grant2014cvx}, CVXPY
\cite{diamond2016cvxpy, agrawal2018rewriting}, CVXR \cite{fu2017cvxr}, and
Convex.jl \cite{udell2014convexjl}.

\paragraph{Solution maps and implicit functions.} An optimization problem can
be viewed as a (possibly multi-valued) function mapping
the problem data, \ie, the numerical data
defining the problem, to the (primal and dual) solution.
This \emph{solution map} is in general
set-valued. In neighborhoods where the solution map is
single-valued, it is an implicit function of the
problem data. In these neighborhoods it becomes meaningful and interesting
to discuss how perturbations in the problem data affect the solution. The point
of this paper is to calculate the effects of such perturbations for
cone programs, in an efficient way.

The use of implicit differentiation to study the sensitivity of solution
mappings of optimization problems dates back several decades, with the works of
Fiacco \cite{fiacco1968} and Robinson \cite{robinson1980} marking significant
milestones. The book \cite{bonnans2000} provides a thorough treatment of the
subject, and \cite{DonRock09} is a good reference for implicit functions more
generally. In recent applications, the framework of implicit functions has been
used to differentiate through quadratic programs \cite{amos2017optnet},
stochastic optimization \cite{donti2017}, physics simulators \cite{de2018},
control algorithms \cite{amos2018differentiable}, and games \cite{ling2018game}.

\paragraph{Automatic differentiation (AD).} Contemporary interest in applications of
implicit differentiation is partly due to the availability of high-quality,
modern, open-source AD software. AD
is a family of techniques for algorithmically computing
exact derivatives of compositions of differentiable functions. Techniques for
AD have been known since at least the 1950s
\cite{beda1959}; see also \cite{wengert1964simple, griewank1989, griewank2008}.
In fact, AD is essentially an efficient way of computing
the chain rule, which can be traced back to a centuries-old manuscript by
Leibniz \cite{leibniz1676} (see \cite{rodriguez2010} for a detailed history).

There are two main variants of AD. Reverse-mode
AD computes the derivative of a composition of atomic
differentiable functions by computing the sensitivity of an output with respect
to the intermediate variables (without materializing the matrices for the
intermediate derivatives). In this way, reverse-mode can efficiently
compute the derivatives of scalar-valued functions. Forward-mode AD
computes the derivative by calculating the sensitivity of the
intermediate variables with respect to an input variable \cite[\S2]{griewank2008}.

Reverse-mode AD was implemented as early as the 1980s
\cite{speelpenning1980compiling}. Its rediscovery by the machine learning
community (where it is known as backpropagation \cite{rumelhart1988}) and the
modern popularity of deep neural networks have led to the development of many
software libraries for reverse-mode AD\@. Examples
include TensorFlow \cite{abadi2016tensorflow, tfe2019}, PyTorch
\cite{paszke2017pytorch}, autograd \cite{maclaurin2015autograd}, and Zygote
\cite{innes2019zygote}. The library JAX \cite{frostig18jax} supports both
reverse-mode and forward-mode AD\@. This requires
representing the derivative of each atomic function as an \emph{abstract linear
map}, \ie, as methods that apply the derivative and its adjoint
\cite{diamond2016matrix}; implementing just
reverse-mode AD only requires the adjoint. Many of the
atomic functions included in these libraries are not differentiable at all
points in their domains. At non-differentiable points, libraries compute
heuristic quantities instead of derivatives. For a discussion on
non-differentiability as it relates to AD, see
\cite[\S14]{griewank2008}.

\paragraph{This paper.} In this paper, we give conditions that guarantee the
existence of the derivative of the solution map for a cone program, and
we provide an expression for the derivative at points where these conditions
are satisfied. As in \cite{bmb18}, our formulation involves expressing the cone
program as the problem of finding a zero of a particular function, specifically
the residual map for a homogeneous self-dual embedding of the program
\cite{ye1994nl, ocpb:16}. We also show how to efficiently compute the
derivative and its adjoint, which involves computing projections onto convex
cones, solving a linear system, and exploiting sparsity.
In \S\ref{sec:implementation}, we present an open-source Python package that
furnishes the derivative of a cone program as an abstract linear map.

\section{The solution map and its derivative}
\label{sec:the_solution_map_and_its_derivative}

We consider a (convex) conic optimization problem
in its primal (P) and dual (D) forms
(see, \eg, \cite[\S4.6.1]{boyd2004convex} or \cite[\S1.4]{ben2001lectures}):
\begin{center}
\begin{tabular}{p{6.5cm}p{6.5cm}}
	{\begin{equation}
			\begin{array}{lll}
				\text{(P)}
				&\minimize  &c\tran x\\
				&\text{subject to} &  Ax+s=b\\
				&  &s\in  \mathcal{K},
			\end{array}
	\end{equation}}
	&
	{ \begin{equation}
			\label{D}
			\begin{array}{lll}
				\text{(D)}&\minimize& b\tran y\\
				&\text{subject to}& A\tran y+c=0\\
				&&y\in  \mathcal{K}^*.
			\end{array}
	\end{equation}}
\end{tabular}
\end{center}
Here
$x\in \RR^n$ is the \emph{primal} variable,
$y\in \RR^m$ is the \emph{dual} variable,
and
$s\in \RR^m$ is the primal \emph{slack} variable.
The set $\mathcal{K}\subseteq \RR^m$
is a nonempty, closed, convex cone with \emph{dual cone}
$\mathcal{K}^*\subseteq \RR^m$.
The
\emph{problem data} are $A \in \RR^{m \times n}$, $b\in \RR^m $, $c\in \RR^n$,
and the cone $ \mathcal{K}$.  (In the sequel, however, we will consider the cone
as fixed.)
In the following we let $N=m+n+1$,
and use $\Pi$ to denote the projection
onto
$\reals^n \times \mathcal{K}^*\times \reals_+$.
Finally, we define
\BEQ
\mathcal{Q}
=\Bigg\{Q = \begin{bmatrix}
0 & A{\tran} & c \\
-A & 0 & b \\
-c{\tran} & -b{\tran} & 0
\end{bmatrix}\in \reals^{N\times N}~\Bigg|~
(A,b,c) \in \reals^{m\times n}\times
\reals^m\times \reals^n \Bigg\}.
\EEQ
Evidently $\mathcal{Q}$ is a proper subspace of the
space of $N\times N$ skew symmetric matrices.
\paragraph{The solution map.}
We call  $(x,y,s)$
a solution of the primal-dual
conic program \eqref{D}
if
\BEQ
\label{eq-kkt}
Ax +s =b, \quad
A\tran y + c = 0, \quad
s \in \mathcal{K}, \quad
y \in \mathcal{K}^*, \quad
s\tran y = 0.
\EEQ
For given problem data,
the corresponding primal-dual conic program \eqref{D} may have
no solution, a unique solution, or multiple solutions.
We focus on the case when it has a unique solution.
We define the \emph{solution map}
$\mathcal{S}:
\reals^{m\times n} \times
\reals^m\times \reals^n \to \reals^{n+2m}$
as the function mapping $(A,b,c)$ to vectors
$(x,y,s)$ that satisfy \eqref{eq-kkt}.
We express $\mathcal S$ as the composition
$\phi \circ s \circ Q$, where
\begin{itemize}
	\item $Q: \reals^{m \times n} \times \reals^m \times \reals^n \to \mathcal Q$
		maps the problem data to the corresponding skew-symmetric matrix in $\mathcal{Q}$,
	\item $s: \mathcal Q \to \reals^N$ furnishes a solution of the homogeneous
		self-dual embedding, and
	\item $\phi: \reals^N \to \reals^n \times \reals^m \times \reals^m$ maps
		a solution of the homogeneous self-dual embedding to the solution
		of the primal-dual pair \eqref{D}.
\end{itemize}
At a point $(A, b, c)$ where $\mathcal S$ is differentiable, the derivative of the
solution map is
\[
\mathsf{D}\mathcal{S}(A,b,c) = \mathsf{D}\phi(z) \mathsf{D}s(Q) \mathsf{D}Q(A,b,c),
\]
by the chain rule. In the remainder of this section, we describe the functions
$Q$, $s$, and $\phi$ and their derivatives, along with sufficient
conditions for their differentiability.

\paragraph{Skew-symmetric mapping.}
Define
\[
	Q = Q(A,b,c)= \begin{bmatrix}
		0 & A{\tran} & c\\
		-A & 0 & b \\
		-c{\tran} & -b{\tran} & 0
	\end{bmatrix}\in \mathcal{Q}.
\]

\paragraph{Homogeneous self-dual embedding.}
The homogeneous self-dual
embedding of \eqref{D}
uses the variable $z\in \reals^N$.
We partition $z$ as $(u,v,w) \in \reals^n \times \reals^m \times \reals$.
The \emph{normalized residual map} introduced in \cite{bmb18} is
the function
$ \mathcal{N}: \RR^N
\times \mathcal{Q}
\to \RR^N$,
defined by
\BEQ
\label{eq-residual}
\mathcal{N}(z,Q)
=\big((Q-\Id)\Pi+\Id\big)(z/|w|).
\EEQ
For a given $Q$, $z$ can be used to
construct the solution
of the primal-dual pair \eqref{D} if and only
if
$\mathcal{N}(z,Q)=0$ and $w>0$ \cite{bmb18}.

The normalized residual map is an affine function
of $Q$,
hence its derivative $\mathsf{D}_Q\mathcal{N}$
always exists, and is given by
\BEQ
\label{e:181203:a}
\mathsf{D}_Q\mathcal{N}(z,Q)(U)=U \Pi(z/|w|),
\quad
\mathsf{D}_Q\mathcal{N}(z,Q)^T(y)=y (\Pi(z/|w|))^T.
\EEQ
We now turn to
$\mathsf{D}_z\mathcal{N}(z,Q)$.
$\mathcal{N}$ is differentiable at $z$,
with $w\neq 0$,
whenever $\Pi$ is differentiable
at $z$, in which case one
 can directly verify that
\BEQ
\label{eq:N:general:deriv}
\mathsf{D}_z\mathcal{N}(z,Q)
=((Q-\Id)\mathsf{D}\Pi(z)+\Id)/w-\sign(w)((Q-\Id)\Pi+\Id)(z/w^2)e^T,
\EEQ
where $e=(0,0,\ldots, 1)\in \reals^N$.  In particular,
when $z$ is a solution of the primal-dual pair \eqref{D},
the second term on the right hand side of 
\eqref{eq:N:general:deriv} vanishes and
we have
\BEQ
\label{e:181203:b}
\mathsf{D}_z\mathcal{N}(z,Q)
=((Q-\Id)\mathsf{D}\Pi(z)+\Id)/w.
\EEQ

\paragraph{Implicit function theorem.}
Suppose that $z$ is a solution
of the primal-dual pair
\eqref{D} for
a given $Q$,
and that $\Pi$ is differentiable at $z$.
Then $\mathcal{N}$
is differentiable at $z$,
$\mathcal{N}(z,Q)=0$
and $w>0$.
Now, suppose that
$\mathsf{D}_z\mathcal{N}(z,Q)$
is invertible.
It follows from the
implicit function theorem (see, \eg, \cite{Dini} and \cite{DonRock09})
that there exists a neighborhood $V\subseteq \mathcal{Q}$
of $Q$ on which the solution $z = s(Q)$ of $\mathcal{N}(z,Q)=0$
is unique. Furthermore, $s$ is differentiable on $V$, $\mathcal{N}(s(Q),Q)=0$
for all
$Q\in V$, and
\BEQ
\label{e:181203:c}
\mathsf{D}s(Q)
=-(\mathsf{D}_z\mathcal{N}(s(Q),Q))^{-1}
\mathsf{D}_Q\mathcal{N}(s(Q),Q).
\EEQ

\paragraph{Solution construction.}
The function $\phi : \reals^{N} \to \reals^{n + 2m}$, given by
\[
\phi(z) = (u, \Pi_\mathcal{K^*}(v),
\Pi_{\mathcal{K}^*}(v)-v)/w,
\]
constructs a solution $(x, y, s)$ of the primal-dual pair \eqref{D} from a
solution $z = (u, v, w)$ of the homogeneous self-dual embedding
($\Pi_{\mathcal{K}^*}$ denotes the projection onto $\mathcal{K}^*$). If
$\Pi_{\mathcal{K}^\star}$ is differentiable with derivative
$\mathsf{D}\Pi_{\mathcal{K}^\star}(v) \in \reals^{m \times m}$, then $\phi$ is
also differentiable, with derivative
\[
\mathsf{D} \phi (z) =
\begin{bmatrix}
I & 0 & -x \\
0 & \mathsf{D}\Pi_{\mathcal{K}^*}(v) & -y \\
0 & \mathsf{D}\Pi_{\mathcal{K}^*}(v) - I & -s,
\end{bmatrix}.
\]
Note that a solution of the homogeneous self-dual embedding can also be
constructed from a solution of the primal-dual pair \cite{bmb18}.

\section{Implementation}
\label{sec:implementation}

In this section we detail how to form
the derivative of a cone program as an abstract linear map.
We also describe our Python package that implements these methods,
and as an example use it to differentiate a semidefinite program.

\paragraph{Sparsity.}
The matrix $A$ is often stored and manipulated as a sparse matrix.
We assume that the sparsity pattern of $A$ is fixed,
meaning we only have to consider the nonzero entries of $A$
when computing the derivative and its adjoint;
this can provide significant speed-ups if $A$ is very sparse.
Of course, one can still furnish the derivative with respect to every
entry of $A$ by making $A$ dense.

\paragraph{Computing the derivative.}
Applying the derivative $\mathsf{D}\mathcal{S}(A, b, c)$ to a
perturbation $(\mathsf{d}A,\mathsf{d}b,\mathsf{d}c)$ corresponds to evaluating
\[
(\mathsf{d}x,\mathsf{d}y,\mathsf{d}s) =
\mathsf{D}\mathcal{S}(A,b,c)(\mathsf{d}A,\mathsf{d}b,\mathsf{d}c)
= \mathsf{D}\phi(z) \mathsf{D}s(Q) \mathsf{D}Q(A,b,c)
(\mathsf{d}A,\mathsf{d}b,\mathsf{d}c),
\]
where $\mathsf{d}A$ has the same sparsity pattern as $A$.
We work from right to left.
The first step is to form
\[
\mathsf{d}Q = \begin{bmatrix}
0 & \mathsf{d}A^T & \mathsf{d}c \\
-\mathsf{d}A & 0 & \mathsf{d}b \\
-\mathsf{d}c^T & -\mathsf{d}b^T & 0 \\
\end{bmatrix}.
\]
The next step is to compute
\[
g = \mathsf{D}_Q\mathcal{N}(s(Q),Q)(\mathsf{d}Q) = \mathsf{d}Q\Pi(z/|w|),
\]
and then
\[
\mathsf{d}z = -M^{-1} g,
\]
where $M=((Q-I)\mathsf{D}\Pi(z)+I)/w$.
One option is to form $M$ as a dense matrix, factorize it,
and then back-solve.
However, when $M$ is large, as is the case in many applications,
this can be impractical.
Instead, we suggest using LSQR \cite{paige1982lsqr} to solve
\BEQ
\begin{array}{ll}
\underset{\mathsf{d}z}{\mbox{minimize}} & \|M\mathsf{d}z + g\|_2^2,
\end{array}
\EEQ
which only requires multiplication by $M$ and $M^T$, and
is of particular interest for semidefinite cone programs.
Next we partition $\mathsf{d}z$ as
$\mathsf{d}z=(\mathsf{d}u,\mathsf{d}v,\mathsf{d}w)$.
The final step is to compute
\[
\begin{bmatrix}
\mathsf{d}x \\
\mathsf{d}y \\
\mathsf{d}s
\end{bmatrix}
=
\begin{bmatrix}
\mathsf{d}u - (\mathsf{d}w) x \\
\mathsf{D}\Pi_{\mathcal{K}^*}(v) \mathsf{d}v - (\mathsf{d}w) y \\
\mathsf{D}\Pi_{\mathcal{K}^*}(v) \mathsf{d}v - \mathsf{d}v - (\mathsf{d}w) s
\end{bmatrix}.
\]

\paragraph{Computing the adjoint of the derivative.}
Applying the adjoint of the derivative to a
perturbation $(\mathsf{d}x,\mathsf{d}y,\mathsf{d}s)$ corresponds to evaluating
\[
(\mathsf{d}A,\mathsf{d}b,\mathsf{d}c) =
\mathsf{D}\mathcal{S}(A,b,c)^T(\mathsf{d}x,\mathsf{d}y,\mathsf{d}s)
= \mathsf{D}Q(A,b,c)^T \mathsf{D}s(Q)^T \mathsf{D}\phi(z)^T
(\mathsf{d}x,\mathsf{d}y,\mathsf{d}s).
\]
We again work right to left.
The first step is to form
\[
\mathsf{d}z
=
\mathsf{D}\phi(z)^T
(\mathsf{d}x,\mathsf{d}y,\mathsf{d}s)
=
\begin{bmatrix}
\mathsf{d}x \\
\mathsf{D}\Pi_{\mathcal{K}^*}^T(v)(\mathsf{d}y + \mathsf{d}s) - \mathsf{d}s \\
-x^T\mathsf{d}x - y^T\mathsf{d}y - s^T\mathsf{d}s
\end{bmatrix}.
\]
Next we form $g = -M^{-T} \mathsf{d}z$, again using LSQR.
Then $\mathsf{d} Q$ is given by
\[
\mathsf{d}Q = g (\Pi(z/|w|))^T.
\]
Instead of explicitly forming $\mathsf{d}Q$, we only obtain
its nonzero entries.
Let its nonzero entries be indexed by $\Omega$;
we compute
\[
(\mathsf{d}Q)_{ij} = \begin{cases}
g_i \Pi(z/|w|)_j & (i,j) \in \Omega \\
0 & \text{otherwise}.
\end{cases}
\]
Partitioning $\mathsf{d} Q$ as
\[
\mathsf{d}Q = \begin{bmatrix}
\mathsf{d}Q_{11} & \mathsf{d}Q_{12} & \mathsf{d}Q_{13} \\
\mathsf{d}Q_{21} & \mathsf{d}Q_{22} & \mathsf{d}Q_{23} \\
\mathsf{d}Q_{31} & \mathsf{d}Q_{32} & \mathsf{d}Q_{33} \\
\end{bmatrix},
\]
the final expressions are given by
\BEAS
\mathsf{d}A &=& \mathsf{d}Q_{12}^T - \mathsf{d}Q_{21} \\
\mathsf{d}b &=& \mathsf{d}Q_{23} - \mathsf{d}Q_{32}^T \\
\mathsf{d}c &=& \mathsf{d}Q_{13} - \mathsf{d}Q_{31}^T.
\EEAS

\paragraph{Integration into AD.}
The calculations that we have described can be immediately
be integrated into the forward and reverse-mode AD software
described in \S\ref{sec:introduction}.
In a forward-mode AD system, one
calculates the sensitivity of the
intermediate variables with respect to perturbations;
this operation corresponds to applying the derivative to those
perturbations.
In a reverse-mode AD system, one computes the sensitivity
of a scalar function with respect
to intermediate variables;
this operation is given by applying the adjoint of the derivative
to the derivative of the scalar function with respect to $x$, $y$, and $s$.

\subsection{Python implementation}
\label{sec:Python_implementation}

We provide a Python implementation of the ideas described in the paper,
available at
\begin{center}
\url{https://www.github.com/cvxgrp/diffcp}.
\end{center}
We use the libraries NumPy for dense linear algebra
\cite{van2011numpy} and
SciPy for sparse linear algebra and its LSQR implementation
\cite{jones2001scipy}.
To solve the cone program, we use the numerical solver SCS \cite{ocpb:16}.
Our implementation supports any cone that
can be expressed as the Cartesian product of the
zero cone, positive orthant, second-order cone,
positive semidefinite cone, exponential cone, and dual exponential cone.
Expressions for the derivative of the projection onto
each of these cones are given in \cite{ali2017semismooth} and \cite{bmb18};
most of these have analytical expressions.

The Python package exposes one function,
\begin{center}
\verb|solve_and_derivative(A, b, c, cone_dict)|,
\end{center}
where \verb|A| is a SciPy sparse matrix,
\verb|b| and \verb|c| are NumPy arrays,
and \verb|cone_dict| is a dictionary representing the cone
$\mathcal K$.
This function returns a solution $(x, y, s)$ of the primal-dual pair and two functions,
\verb|derivative(dA, db, dc)|
and
\verb|adjoint_derivative(dx, dy, ds)|,
which respectively apply the derivative and its adjoint (at \verb|(A, b, c)|) to
their inputs and return the result.

\subsection{Example}
\label{sec:differentiating_a_semidefinite_program}

We consider differentiating a semidefinite program
\BEQ
\begin{array}{ll}
\mbox{minimize} & \mathbf{tr}(C^TX) \\
\mbox{subject to} & X \succeq 0, \\
& \mathbf{tr}(A_i X) = b_i, \quad i=i,\ldots,p,
\end{array}
\label{eq:sdp}
\EEQ
with optimization variable $X\in\symm_+^n$ and problem data
$C\in\symm^n$, $A_1,\ldots,A_p\in\symm^n$, and
$b\in\reals^p$ ($\symm^n$ denotes the set of symmetric matrices in $\reals^{n
\times n}$, and $\symm_+^n$ denotes the set of symmetric positive semidefinite
matrices in $\reals^{n \times n}$). This problem can be readily cast
as a standard cone program \eqref{D},
where the cone is the Cartesian product of a zero cone
and a semidefinite cone.

We generated a feasible, bounded
random instance of \eqref{eq:sdp} with $p=100$ and $n=300$.
We used our Python package to solve this instance and retrieve the derivative
and its adjoint. The solve, which calls into SCS, took about 195.7 seconds on a
machine with a six-core Intel i7-8700K. Next, we computed the derivative of
the optimal value of \eqref{eq:sdp} with respect to each of the $A_i$ by
applying the adjoint of the derivative to the gradient of the objective.
Applying the adjoint of the derivative took about 191.7 seconds, which is
roughly equal to the time it took to solve the problem. Note that
we calculated the derivative of the objective with respect to 4,515,000
elements (each of the $A_i$), which required the solution of a 90,401 $\times$
90,401 linear system.
This evidently would not have been practical had we not treated the
linear mapping as an abstract operator.

\section*{Acknowledgements}
We thank Brandon Amos and Zico Kolter, who concurrently and independently
derived calculations similar to the ones in this paper, for many useful
discussions. Shane Barratt is supported by the National Science Foundation
Graduate Research Fellowship under Grant No. DGE-1656518.
Walaa Moursi is partially supported by the 
Natural Science and 
Engineering Research Council of Canada 
Postdoctoral Fellowship. 
\bibliography{conic_derivative}

\end{document}

%% file: defs.tex
\newcommand{\BEAS}{\begin{eqnarray*}}
\newcommand{\EEAS}{\end{eqnarray*}}
\newcommand{\BEQ}{\begin{equation}}
\newcommand{\EEQ}{\end{equation}}
\newcommand{\BIT}{\begin{itemize}}
\newcommand{\EIT}{\end{itemize}}

\newcommand{\eg}{{\it e.g.}}
\newcommand{\ie}{{\it i.e.}}

\newcommand{\reals}{{\mbox{\bf R}}}
\newcommand{\symm}{{\mbox{\bf S}}}

{\begin{quote}}{\end{quote}}

\makeatletter
\long\def\@makecaption#1#2{
   \vskip 9pt
   \begin{small}
   \setbox\@tempboxa\hbox{{\bf #1:} #2}
   \ifdim \wd\@tempboxa > 5.5in
        \begin{center}
        \begin{minipage}[t]{5.5in}
        \addtolength{\baselineskip}{-0.95pt}
        {\bf #1:} #2 \par
        \addtolength{\baselineskip}{0.95pt}
        \end{minipage}
        \end{center}
   \else
	\hbox to\hsize{\hfil\box\@tempboxa\hfil}
   \fi
   \end{small}\par
}
\makeatother

\newcounter{oursection}

\newcounter{lecture}